\documentclass[12pt]{amsart}

\usepackage{amsmath,amssymb,graphicx,mathrsfs,overpic,tikz-cd}

\newtheorem{lemma}{Lemma}[section]
\newtheorem{proposition}[lemma]{Proposition}
\newtheorem{theorem}[lemma]{Theorem}

\newtheorem{conjecture}[lemma]{Conjecture}

\theoremstyle{plain}

\theoremstyle{plain}

\def\C{\mathbb C}
\def\R{\mathbb{R}}
\def\Z{\mathbb{Z}}
\def\to{\rightarrow}
\def\SL{\mathrm{SL}}
\def\tr{\mathrm{tr}}
\def\S{\mathbb{S}}

\def\A{\mathcal{A}}
\def\I{\mathcal{I}}
\def\l{\langle\!\langle}
\def\r{\rangle\!\rangle}

\def\g{\mathbf{g}}
\def\i{\mathbf{i}}

\title{Trace-free characters and abelian knot contact homology II}

\author{Fumikazu Nagasato and Shinnosuke Suzuki}

\address{Department of Mathematics, Meijo University, 
Tempaku, Nagoya 468-8502, Japan}
\email{fukky@meijo-u.ac.jp}

\address{Toyota Technical Development Corporation (TTDC), 
1-9 Imae, Hanamoto-cho, Toyota, Aichi 470-0334, Japan}
\email{}

\subjclass[2020]{Primary 57K31; Secondary 57K18}

\keywords{abelian knot contact homology, character varieties, ghost characters, trace-free characters}

\begin{document}

\begin{abstract}
We show that the $(4,5)$- and $(5,6)$-torus knots admit ghost characters. 
Consequently, these knots provide counterexamples to Ng's conjecture, 
which proposes an isomorphism between the complexification of degree $0$ 
abelian knot contact homology and the coordinate ring of the character variety 
of the $2$-fold branched cover of the $3$-sphere branched along a knot. 
While Ng's conjecture has been verified for all $2$-bridge and $3$-bridge knots, 
we demonstrate, via ghost characters, how this isomorphism fails for these torus knots. 
\end{abstract}
\maketitle


\section{Introduction}\label{intro}
In \cite{Ng1}, L. Ng introduced degree $0$ knot contact homology $HC_0(K)$, 
which is a noncommutative algebra over $\mathbb{Z}$ defined for a knot $K$. 
He further focused on its abelianization, 
known as degree $0$ abelian knot contact homology $HC^{ab}_0(K)$,  
and conjectured that the complexification of $HC^{ab}_0(K)$ is isomorphic to 
the coordinate ring $\mathbf{C}[X(\Sigma_2 K)]$ of the character variety 
of the $2$-fold branched cover $\Sigma_2 K$ of the $3$-sphere 
$\mathbb{S}^3$ branched along $K$ (see Conjecture \ref{conj_ng}).
As shown in \cite{Nagasato5}, the structure underlying this correspondence 
can be described in terms of the trace-free slice $S_0(K)$ 
and the fundamental variety $F_2(K)$, together with the maps 
\[
\begin{tikzcd}
S_0(K) \arrow[r, "\widehat{\Phi}"] \arrow[d, "q"'] & X(\Sigma_2 K) \arrow[ld, "h^*"] \\
F_2(K) &
\end{tikzcd}
\]
In this setting, we introduced in \cite[Definition 4.7]{Nagasato5} 
a complete obstruction to Ng's conjecture, called a ghost character of a knot. 
As a consequence, Ng's conjecture holds for a knot $K$ 
if and only if $K$ admits no ghost characters \cite[Theorem 4.10]{Nagasato5}. 
In the present paper, we show that the $(4,5)$- and $(5,6)$-torus knots 
admit ghost characters and hence provide counterexamples to Ng's conjecture. 

This paper is organized as follows. 
Section \ref{overview} gives an overview of Ng's conjecture via the maps 
$h^*$, $\widehat{\Phi}$, and $q$, together with brief reviews of degree 0 
abelian knot contact homology, the character variety of the 2-fold branched cover 
of the $3$-sphere branched along a knot, the trace-free slice, the fundamental variety, 
and ghost characters. 
In Section \ref{ghost}, we determine the ghost characters of the $(4,5)$- and 
$(5,6)$-torus knots, thereby showing that Ng's conjecture does not hold for these knots. 
In Section \ref{sec_h}, we investigate in detail how the isomorphic correspondence 
proposed in Ng's conjecture collapses for these knots through the maps $h^*$ and $\widehat{\Phi}$. 


\section{Overview of Ng's conjecture}\label{overview}
Before turning to the details, we briefly review the general framework of Ng's conjecture. 
For further details, see \cite{Nagasato5}. 

Let $K \subset \R^3$ be a knot with an $n$-crossing diagram $D_K$, and let 
\[
G(K)=\langle m_1,\cdots,m_n \mid r_1, \cdots, r_{n-1} \rangle
\] 
be the associated Wirtinger presentation of the knot group. 
Let $\A_n^{ab}$ denote the polynomial ring over $\Z$ 
generated by indeterminates $a_{ij}$ $(1 \le i < j \le n)$, 
with the convention that $a_{ii}=-2$ for all $i$, where the indices correspond to the meridians 
$m_1,\cdots,m_n$. 
We denote by $\mathcal{I}_{D_K}$ the ideal\footnote{Although the original definition of the ideal 
$\I_{D_K}$ is given in a different form, it is equivalent to the ideal defined in the present paper. 
For details, see \cite[Section 6]{Nagasato3}.}  
of $\mathcal{A}_n^{ab}$ generated by the elements 
\[
a_{lj} + a_{lk} + a_{li} a_{ij}, 
\]
where $(i,j,k)$ corresponds to a crossing of $D_K$ such that $i$ is the overarc  
and $j,k$ are the underarcs, and $1 \leq l \leq n$. 
Such a triple $(i,j,k)$ with $i < k$ is called a Wirtinger triple of $D_K$. Explicitly, 
\[
\mathcal{I}_{D_K}
:=
\left\langle
a_{lj} + a_{lk} + a_{li} a_{ij}
\ \middle|\
\begin{array}{l}
(i,j,k): \text{any Wirtinger triple of } D_K, \\
1 \leq l \leq n
\end{array}
\right\rangle.
\]
With this notation, the degree 0 abelian knot contact homology of $K$,
denoted by $HC^{ab}_0(K)$, is defined as the quotient
\[
HC^{ab}_0(K)
:=
\mathcal{A}_n^{ab} / \mathcal{I}_{D_K}.
\]
It is known that $HC^{ab}_0(K)$ is a knot invariant, up to ring isomorphism fixing
$\mathbb{Z}$ pointwise (see \cite{Ng1,Ng2}).

Let $K$ be a knot in the $3$-sphere $\mathbb{S}^3$. 
Suppose that the knot group $G(K)$ admits a presentation\footnote{This presentation is not necessarily 
a Wirtinger presentation.}
\[
G(K)=\langle m_1,\dots,m_n \mid r_1,\dots,r_{n-1} \rangle,
\]
where $m_1,\dots,m_n$ are meridians of $K$.
Let $\mathfrak{S}_0(K)$ denote the set of the characters $\chi_\rho$ of 
trace-free representations 
$\rho: G(K) \to \SL_2(\mathbb{C})$, that is, representations satisfying
$\tr(\rho(\mu_K)) = 0$ for a meridian $\mu_K$ of $K$. 
Following \cite{Culler-Shalen,Nagasato5}, 
the set $\mathfrak{S}_0(K)$ is embedded as a closed algebraic set in 
$\mathbb{C}^{{n \choose 2}+{n \choose 3}}$ via the trace map
\[
\tilde{t} : \mathfrak{S}_0(K) \rightarrow \mathbb{C}^{{n \choose 2}+{n \choose 3}},
\ 
\tilde{t}(\chi_\rho)
:=
\left(x_{ij}(\chi_\rho);\, x_{ijk}(\chi_\rho)\right),
\]
where the parameters $x_{ij}$ and $x_{ijk}$ are defined using negative traces by
\[
\begin{array}{ll}
x_{ij}(\chi_\rho) := - t_{m_i m_j}(\rho) =-\tr(\rho(m_i m_j))  & (1 \le i < j \le n), \\
x_{ijk}(\chi_\rho) := - t_{m_i m_j m_k}(\rho)=-tr(\rho(m_i m_j m_k)) & (1 \le i < j < k \le n).
\end{array}
\]
We call the closed algebraic set $\tilde{t}(\mathfrak{S}_0(K))$ the trace-free slice 
of $K$, and denote it by $S_0(K)$. The parametrization of $S_0(K)$ 
depends on the choice of a generating set of $G(K)$, but only up to biregular equivalence. 
Therefore, $S_0(K)$ is an invariant of a knot $K$ up to biregular equivalence. 

With the above notation, let
\[
p \colon C_2 K \longrightarrow E_K
\]
be the $2$-fold cyclic cover of the knot exterior $E_K$, chosen so that the image 
$p(\mu_2)$ of a meridian $\mu_2$ of $C_2 K$ is homotopic to $m_1^2$. 
The $2$-fold branched cover $\Sigma_2 K$ of $\mathbb{S}^3$ branched along $K$ is obtained  
from $C_2K$ by trivially filling a solid torus along the boundary of $C_2K$, 
that is, by gluing the standard meridian of the solid torus to $\mu_2$. 
The covering map $p$ induces a natural injection 
\[
p_* \colon \pi_1(C_2 K) \longrightarrow G(K), 
\]
and hence there exists an isomorphism 
\[
\pi_1(\Sigma_2 K) \cong \mathrm{Im}(p_*) / \l m_1^2 \r 
\]
where $\l * \r$ denotes the normal closure of the group $\langle * \rangle$. 
In particular, by a result of Fox \cite{Fox}\footnote{Fox's theorem applies to presentations of knot 
groups generated by meridians.} 
(see also \cite{Kinoshita} and \cite[Theorem 4.5]{Nagasato5}),
which is based on a coset decomposition of $\mathrm{Im}(p_*)$ using the Schreier system
$\{1,m_1\}$, we obtain the presentation
\[
\pi_1(\Sigma_2K) \cong \langle m_1m_i\ (2 \leq i \leq n) \mid w(r_j),w(m_1r_jm_1^{-1})\ 
(1 \leq j \leq n-1), m_i^2\ (1 \leq i \leq n) \rangle,
\] 
where $w(r_j)$ (resp. $w(m_1r_jm_1^{-1})$) denotes the word obtained by rewriting 
$r_j$ (resp. $m_1r_jm_1^{-1}$) in terms of the generators $m_1m_2, \cdots, m_1m_n$. 

It is shown in \cite{Culler-Shalen, Fricke, Horowitz, Vogt} 
(see also \cite{Gonzalez-Montesinos}) that, for any element 
$g \in \pi_1(\Sigma_2 K)$ and any representation 
$\rho_* : \pi_1(\Sigma_2 K) \rightarrow \SL_2(\mathbb{C})$, 
the trace $\tr(\rho_*(g))$ can be written as a polynomial in the traces
\[
\begin{array}{ll}
z_{ab}(\chi_{\rho_*}):=t_{(m_am_b)}(\chi_{\rho_*}) &(1 \leq a < b \leq n),\\ 
y_{def}(\chi_{\rho_*}):=t_{(m_1m_d)(m_1m_e)(m_1m_f)}(\chi_{\rho_*}) & (2 \leq d < e < f \leq n),  
\end{array}
\]
where $t_g(\chi_{\rho}):=\tr(\rho(g))$ denotes the trace function. 
Again, following \cite{Culler-Shalen,Nagasato5}, the character variety $X(\Sigma_2 K)$ is defined 
as the image of the trace map 
\[
t:\mathfrak{X}(\Sigma_2K) \to \C^{{n \choose 2}+{n-1 \choose 3}},
\ t(\chi_{\rho_*}):=(z_{ab}(\chi_{\rho_*});y_{def}(\chi_{\rho_*})). 
\]
Accordingly, $X(\Sigma_2 K)$ admits the parametrization 
\[
X(\Sigma_2K) = 
\left\{
(z_{ab}(\chi_{\rho_*});y_{def}(\chi_{\rho_*})) \in \C^{{n \choose 2}+{n-1 \choose 3}} \mid 
\chi_{\rho_*} \in \mathfrak{X}(\Sigma_2K)
\right\}
\]

Using the above formulation, we recall the map introduced in \cite{Nagasato-Yamaguchi}
\[
\widehat{\Phi} \colon \mathfrak{S}_0(K) \longrightarrow \mathfrak{X}(\Sigma_2 K),
\]
defined as follows. 
For a trace-free character $\chi_\rho \in \mathfrak{S}_0(K)$, 
the image $\widehat{\Phi}(\chi_\rho)$ is the character of the representation 
$\rho_* : \pi_1(\Sigma_2 K) \to \SL_2(\mathbb{C})$ given by 
\begin{eqnarray}\label{rep_s2k}
\rho_*(g):=
(\sqrt{-1})^{\alpha(p_*(g))}\,\rho(p_*(g)),
\ g \in \pi_1(\Sigma_2 K),
\end{eqnarray}
where $\alpha : G(K) \to H_1(E_K)=\langle \mu_K \rangle \cong \mathbb{Z}$ 
denotes the abelianization. 
Under the parametrizations for $S_0(K)$ and $X(\Sigma_2 K)$ introduced above,
the map $\widehat{\Phi}$ is expressed as a polynomial map (see \cite{Nagasato5} for details).
More precisely, for a trace-free character $\chi_\rho = (x_{ij}; x_{ijk}) \in S_0(K)$, we have
\[
\widehat{\Phi}(x_{ij}; x_{ijk})
=
\left(
x_{ab};\,
x_{1d}x_{1e}x_{1f}
-\frac{1}{2}
\bigl(
x_{1d}x_{ef}
+x_{1e}x_{df}
+x_{1f}x_{de}
\bigr)
\right) \in \C^{{n \choose 2}+{n-1 \choose 3}},
\]
where $1 \le a < b \le n$ and $2 \le d < e < f \le n$. This expression defines the polynomial map 
$\widehat{\Phi}: S_0(K) \to X(\Sigma_2K)$. 

The coordinate ring $\mathbf{C}[X(\Sigma_2 K)]$ of the closed algebraic set $X(\Sigma_2 K)$, 
with the coordinates $(z_{ab}; y_{def})$, is defined as the ring of regular functions 
on $X(\Sigma_2 K)$, which is isomorphic to the quotient: 
\[
\mathbf{C}[X(\Sigma_2 K)]
\cong
\C[z_{ab}; y_{def}]/\sqrt{I_{X(\Sigma_2 K)}},
\]
where $I_{X(\Sigma_2 K)}$ denotes the ideal of polynomials in $\C[z_{ab}; y_{def}]$ 
vanishing on $X(\Sigma_2 K)$, and $\sqrt{I_{X(\Sigma_2 K)}}$ is its radical.

With this algebraic setup, we are now ready to state Ng's conjecture precisely.  

\begin{conjecture}[Conjecture 5.7 in \cite{Ng2}]
\label{conj_ng}
Let $G(K)=\langle m_1,\cdots, m_n \mid r_1,\cdots, r_{n-1} \rangle$ 
be a Wirtinger presentation. Then the ring homomorphism
\[
g: HC_0^{ab}(K) \otimes \C \to \mathbf{C}[X(\Sigma_2K)]
\]
given by $g(a_{ij}) = -t_{m_i m_j}$ $(1 \leq i < j \leq n)$, $g(1)=1$, is an isomorphism. 
\end{conjecture}

Since the coordinate ring of any closed algebraic set is reduced, 
it is natural to consider the nilradical quotient of $HC_0^{ab}(K) \otimes \C$ 
in the above conjecture. 
In \cite[Theorem 4.9 (1)]{Nagasato5}, it was shown that Ng's conjecture, 
under this modification, holds for all 2-bridge\footnote{This case was originally shown by Ng \cite{Ng2}.} 
and 3-bridge knots. 
More generally, the modified conjecture holds for a knot $K$ if and only if $K$ admits 
no ghost characters \cite[Theorem 4.10]{Nagasato5}. 

As discussed above, the notion of a ghost character gives a complete obstruction 
to the modified conjecture. 
We now recall the definition of a ghost character for a knot $K$ in $\S^3$. 
In \cite[Theorem 2.1]{Nagasato5}, it was shown that $S_0(K)$,  
associated with a Wirtinger presentation 
\[
G(K)=\langle m_1,\cdots,m_n \mid r_1,\cdots,r_{n-1} \rangle, 
\]
is isomorphic to the common solutions of the following equations: 
\begin{description}
\item[(F2)] the fundamental relations 
\begin{eqnarray*}
&x_{ak}=x_{ij}x_{ai}-x_{aj}&\\
&(1 \leq a \leq n,\ (i,j,k):\mbox{a Wirtinger triple}),&
\end{eqnarray*}
\item[(GH)] the general hexagon relations 
\begin{eqnarray*}
&x_{i_1 i_2 i_3} \cdot x_{j_1 j_2 j_3}
=\frac{1}{2}
\left|
\begin{array}{ccc}
x_{i_1 j_1} & x_{i_1 j_2} & x_{i_1 j_3}\\
x_{i_2 j_1} & x_{i_2 j_2} & x_{i_2 j_3}\\
x_{i_3 j_1} & x_{i_3 j_2} & x_{i_3 j_3}
\end{array}\right|,&\\
&(1 \leq i_1<i_2<i_3 \leq n,\ 1 \leq j_1<j_2<j_3 \leq n),&
\end{eqnarray*}
\end{description}
with the convention $x_{ii}=2$, $x_{ji}=x_{ij}$ and 
$x_{i_{\sigma(1)}i_{\sigma(2)}i_{\sigma(3)}}=\mathrm{sign}(\sigma)x_{i_1i_2i_3}$ 
for any permutation $\sigma \in \mathfrak{S}_3$. 
The fundamental relations (F2) define a closed algebraic set
\[
F_2(K)=\left\{(x_{12},\cdots,x_{n-1n}) \in \C^{n \choose 2}\ \left|\ 
\begin{array}{c}
x_{ak}=x_{ij}x_{ai}-x_{aj} \mbox{ for any $1 \leq a \leq n$}\\
\mbox{and any Wirtinger triple $(i,j,k)$}
\end{array}
\right.\right\}. 
\]
We refer to $F_2(K)$ as the fundamental variety of a knot $K$. 
By forgetting the last $n \choose 3$ coordinates $(x_{ijk})$ of the trace-free slice $S_0(K)$, 
we obtain a natural projection 
\[
q: S_0(K) \to F_2(K) \subset \C^{n \choose 2}.
\] 
In this setting, a ghost character of $K$ is defined to be a point $(x_{ij}) \in F_2(K)$ 
that does not satisfy one of (GH), and hence does not lift to $S_0(K)$. 

The coordinate ring $\mathbf{C}[F_2(K)]$ of the fundamental variety $F_2(K)$
is given by
\[ 
\mathbf{C}[F_2(K)]=\frac{\C[x_{12},\cdots,x_{n n-1}]}
{\sqrt{\langle x_{aj}+x_{ak}-x_{ai}x_{ij}\ \mid \mbox{\small $(i,j,k)$: any Wirtinger triple, } 
\mbox{$1 \leq a \leq n$} \rangle}}. 
\]
Consequently, the map
\[
f : HC_0^{ab}(K) \otimes \C \rightarrow \mathbf{C}[F_2(K)],\ 
f(a_{ij})=-x_{ij}\ (1 \le i < j \le n),\ f(1)=1,
\]
is a well-defined ring homomorphism (see also \cite[Proposition 4.2]{Nagasato5}).
Moreover, the kernel of $f$ is obviously the nilradical $\sqrt{0}$. 
Therefore, for any knot $K$, the quotient ring 
$(HC_0^{ab}(K) \otimes \C)/\sqrt{0}$ is isomorphic to $\mathbf{C}[F_2(K)]$.

Accordingly, Conjecture \ref{conj_ng} in its nilradical quotient claims that, 
for any knot $K$, the ring homomorphism 
\[
h \colon \mathbf{C}[F_2(K)] \longrightarrow \mathbf{C}[X(\Sigma_2 K)],
\]
defined by $h(x_{ij}) = t_{m_i m_j}$ for $1 \le i < j \le n$ and $h(1)=1$, is an isomorphism. 
By Hilbert's Nullstellensatz, this is equivalent to the following statement.
\begin{conjecture}[Conjecture 4.4 in \cite{Nagasato5}]
\label{conj_nag}
For any knot $K$, the pull-back of $h$: 
\[
h^* : X(\Sigma_2 K) \longrightarrow F_2(K),
\]
given by $v(h^*(z)) = h(v)(z)$ for $v \in \mathbf{C}[F_2(K)]$ and $z \in X(\Sigma_2 K)$, 
is an isomorphism of algebraic sets. 
\end{conjecture}
It follows from \cite[Section 4]{Nagasato5} that the map $h^*$ is expressed explicitly 
as a polynomial map by 
\[
h^*((z_{ab}; y_{abc})) = (z_{ab}). 
\]

In this paper, we show that the $(4,5)$-torus knot $T_{4,5}$ and the $(5,6)$-torus knot $T_{5,6}$ 
admit ghost characters, and hence provide counterexamples to Conjecture \ref{conj_nag}, 
by Theorem 4.10 in \cite{Nagasato5}. 
More precisely, we obtain the following result. 

\begin{theorem}[cf. Theorems \ref{thm_main2} and \ref{thm_main}]\label{thm_ng-conj}
The map $h^*$ is surjective but not injective for $T_{4,5}$,  
and neither surjective nor injective for $T_{5,6}$. 
In particular, Conjecture $\ref{conj_nag}$ (Conjecture $\ref{conj_ng}$ 
with the nilradical quotient) does not hold for $T_{4,5}$ and $T_{5,6}$. 
\end{theorem}

The strategy for proving Theorem \ref{thm_ng-conj} is to determine the ghost characters 
for $T_{4,5}$ and $T_{5,6}$. 
Indeed, for any knot $K$, the preimage of a ghost character under the map $h^*$ 
must consist of either two points lying outside ${\rm Im}(\widehat{\Phi})$ or be empty 
(see the proof of Theorem 4.10 in \cite{Nagasato5} for details). 
On the other hand, if a point $\mathbf{x} \in F_2(K)$ is not a ghost character, 
then the preimage $(h^*)^{-1}(\mathbf{x})$ consists of exactly one character in $X(\Sigma_2 K)$, 
which lies in ${\rm Im}(\widehat{\Phi})$ (see the proof of \cite[Theorem 4.9]{Nagasato5}). 


\section{Ghost characters of a knot and Conjecture \ref{conj_nag}}\label{ghost}
As shown in \cite[Theorem 4.8]{Nagasato5}, any knot with bridge index less than 4 
admits no ghost characters. 
Computational evidence indicates that $T_{4,5}$ and $T_{5,6}$, 
which are 4-bridge and 5-bridge knots, respectively, admit ghost characters. 
In this section, we carry out explicit computations of the ghost characters 
for these knots, largely by hand. 


\subsection{Ghost characters of the (4,5)-torus knot}\label{subsec_ghost_t45}
\begin{figure}[htbp]
\[
D=\begin{minipage}{11cm}
\begin{overpic}[width=\hsize]{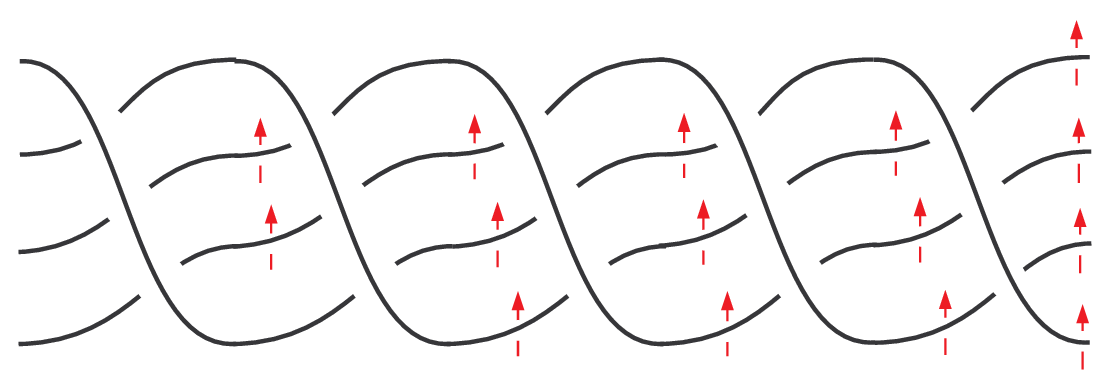}
\put(99,5){$m_1$}
\put(99,14){$m_2$}
\put(99,22.5){$m_3$}
\put(99,31.5){$m_4$}
\put(79,7){$m_5$}
\put(76.5,16){$m_6$}
\put(75,24){$m_7$}
\put(59,7){$m_8$}
\put(57,16){$m_9$}
\put(54.5,24){$m_{10}$}
\put(39,7){$m_{11}$}
\put(37,16){$m_{12}$}
\put(35,24){$m_{13}$}
\put(17,16){$m_{14}$}
\put(16,24){$m_{15}$}
\end{overpic}
\end{minipage}
\]
\caption{A diagram of $T_{4,5}$ in braid (bridge) position and meridians $m_1,\cdots,m_{15}$. 
The four parallel curves connecting the both sides of the diagram are omitted.}
\label{diagram_T45}
\end{figure}

Let $D$ be the diagram of $T_{4,5}$ shown in Figure \ref{diagram_T45}. 
We assign meridians $m_1,\cdots,m_{15}$ as indicated in Figure $\ref{diagram_T45}$, 
thereby obtaining the Wirtinger presentation of $G(T_{4,5})$ associated with $D$. 

With this notation, the fundamental variety $F_2(T_{4,5})$ is given by 
\[
F_2(T_{4,5})=\left\{(x_{12},\cdots,x_{14,15}) \in \C^{15 \choose 2}\ \left|\ 
\begin{array}{c}
x_{ak}=x_{ij}x_{ai}-x_{aj} \mbox{ for any $1 \leq a \leq 15$}\\
\mbox{and any Wirtinger triple $(i,j,k)$}
\end{array}
\right.\right\}. 
\]
As shown in \cite{Nagasato5}, a knot $K$ in bridge position generally admits 
a systematic elimination process for the fundamental relations (F2). 
This method can also be applied to a knot in braid position\footnote{The method for knots 
in braid position was originally developed in \cite{Nagasato3}; 
see also a similar approach in \cite{Ng1,Ng2}.}.  
We demonstrate this process in the present case $T_{4,5}$. 

First, for $1 \leq a \leq 15$, we have the following fundamental relations (F2):
\[
\begin{array}{ll}
x_{a,15}=x_{8,11}x_{a,11}-x_{a,8}, & x_{a,14}=x_{11,13}x_{a,11}-x_{a,13},\\
x_{a,13}=x_{5,8}x_{a,8}-x_{a,5}, & x_{a,12}=x_{8,10}x_{a,8}-x_{a,10},\\
x_{a,11}=x_{8,9}x_{a,8}-x_{a,9}, & x_{a,10}=x_{1,5}x_{a,5}-x_{a,1},\\
x_{a,9}=x_{5,7}x_{a,5}-x_{a,7},  & x_{a,8}=x_{5,6}x_{a,5}-x_{a,6},\\
x_{a,7}=x_{1,4}x_{a,1}-x_{a,4}, & x_{a,6}=x_{1,3}x_{a,1}-x_{a,3},\\
x_{a,5}=x_{1,2}x_{a,1}-x_{a,2}, & \\
x_{a,4}=x_{11,12}x_{a,11}-x_{a,12}, & x_{a,3}=x_{4,11}x_{a,4}-x_{a,11},\\
x_{a,2}=x_{4,15}x_{a,4}-x_{a,15}, & x_{a,1}=x_{4,14}x_{a,4}-x_{a,14}.
\end{array}
\]
Note that the last 4 types of relations (F2) are written for triples $(i,j,k)$ with $j>k$, 
using the symmetry between $j$ and $k$, for technical reasons.  
Although this system can be solved by computer, 
we explain how the elimination can be carried out mostly by hand,  
following the strategy in the proof of \cite[Theorem 4.8]{Nagasato5}. 

To facilitate the elimination process, we introduce the following additional relations
derived from (F2), which will also be referred to as the fundamental relations (F2) below: 
\[
\begin{array}{ll}
\fbox{$x_{12}$} &=x_{4,14}x_{2,4}-x_{2,14}
=x_{4,14}x_{4,15}-(x_{4,15}x_{4,14}-x_{14,15})=\fbox{$x_{14,15}$},\\
\fbox{$x_{13}$} &=x_{4,11}x_{1,4}-x_{1,11}=x_{4,11}x_{4,14}-(x_{4,14}x_{4,11}-x_{11,14})
=\fbox{$x_{11,14}$},\\
\fbox{$x_{14}$} &=x_{4,14}=\fbox{$x_{11,12}x_{11,14}-x_{12,14}$},\\
\fbox{$x_{23}$} &=x_{4,11}x_{2,4}-x_{2,11}=x_{4,11}x_{4,15}-(x_{4,15}x_{4,11}-x_{11,15})
=\fbox{$x_{11,15}$},\\
\fbox{$x_{24}$} &=x_{4,15}=\fbox{$x_{11,12}x_{11,15}-x_{12,15}$},\\
\fbox{$x_{34}$} &=x_{4,11}=\fbox{$x_{11,12}$}.
\end{array}
\]
The basic idea for reducing the parameters of $F_2(K)$ 
is to repeatedly apply the relations (F2) as follows. 
We first eliminate $x_{a,15}$ $(1 \leq a \leq 14)$ from (F2) using 
\[
x_{a,15}=x_{8,11}x_{a,11}-x_{a,8}.  
\]
Next, we eliminate $x_{a,14}$ $(1 \leq a \leq 13)$ from (F2), including the equations 
obtained in the previous step, using the relations: 
\[
x_{a,14}=x_{11,13}x_{a,11}-x_{a,13}. 
\]
Iterating this procedure down to $x_{a,5}$ $(1 \leq a \leq 4)$, 
we eventually express $x_{a,15}$ $(1 \leq a \leq 14)$, $\cdots$, $x_{a,5}$ $(1 \leq a \leq 4)$ 
as elements of the polynomial ring 
\[
R=\C[x_{12},x_{13},x_{14},x_{23},x_{24},x_{34}].
\] 
This recursive elimination substantially simplifies the parametrization of $S_0(T_{4,5})$. 
The resulting parametrization corresponds to the following $4$-bridge knot group 
presentation of $T_{4,5}$ (see the presentation in (\ref{t45-4bridge-2}) for more precise description): 
\begin{eqnarray}\label{t45-4bridge}
G(T_{4,5}) \cong \langle m_1,\ldots,m_4 \mid w_1,w_2,w_3 \rangle.
\end{eqnarray}

In fact, as discussed in the proof of \cite[Theorem 4.8]{Nagasato5}, 
the above elimination process can be understood naturally 
as the following diagrammatic operation, which reflects the diagram $D$ in braid (bridge) position. 
We first interpret $x_{ak}$ as a loop $s_{ak}$ in the knot exterior $E_{T_{4,5}}$ 
freely homotopic to $m_a m_k$, 
and then decompose $s_{ak}$ into two arcs $c_a$ and $c_k$, corresponding to $m_a$ and $m_k$, 
respectively. Next, keeping $c_a$ and the endpoints of $c_k$ fixed, 
we slide $c_k$ across the crossing along the $k$th arc of $D$. 
\[
\begin{minipage}{4cm}
\begin{overpic}[width=\hsize]{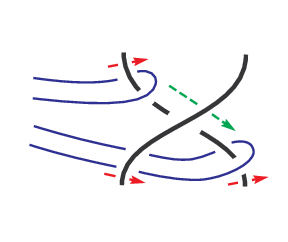}
\put(10,65){$s_{ak}$}
\put(40,10){$i$}
\put(80,10){$j$}
\put(40,70){$k$}
\put(85,48){\small sliding}
{}\end{overpic}
\end{minipage}
=
\begin{minipage}{3.5cm}
\begin{overpic}[width=\hsize]{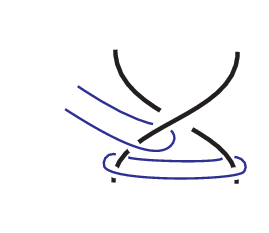}
\put(10,55){$s_{ai}$}
\put(40,8){$i$}
\put(60,10){$s_{ij}$}
\put(85,8){$j$}
\put(38,72){$k$}
\end{overpic}
\end{minipage}
-
\begin{minipage}{3.5cm}
\begin{overpic}[width=\hsize]{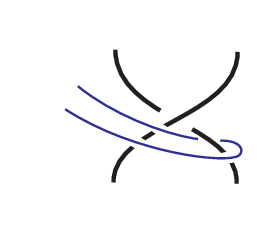}
\put(10,55){$s_{aj}$}
\put(40,8){$i$}
\put(85,8){$j$}
\put(38,72){$k$}
\end{overpic}
\end{minipage}.
\]
(In the case of the current diagram $D$ in Figure \ref{diagram_T45}, slide $c_k$ from left to right.) 
Finally, we resolve the winding part, namely, the part of the resulting loop 
lying under the $i$th arc of $D$, using the trace-free Kauffman bracket skein relation:
\[
\begin{minipage}{10cm}
\begin{overpic}[width=\hsize]{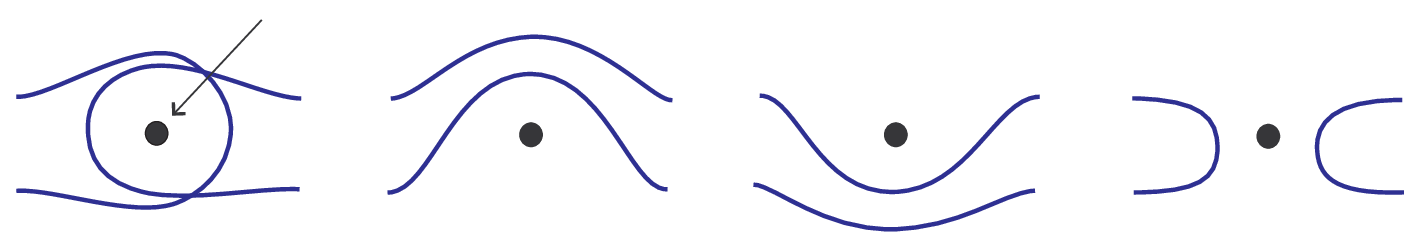}
\put(22.5,7){$=$}
\put(49.5,7){$+$}
\put(75,7){$+$}
\put(19.5,15){\small knot}
\end{overpic}
\end{minipage}.
\]
The resulting loops $s_{ai}$, $s_{aj}$, and $s_{ij}$ are freely homotopic to 
$m_am_i$, $m_am_j$, and $m_im_j$, repsectively. 
Interpreting a disjoint union of loops as the product of the corresponding monomials, 
we recover the relation (F2): $x_{ak}=x_{ij}x_{ai}-x_{aj}$. 
This operation may be considered as an action of braids on the polynomial ring $R$. 
From this viewpoint, we obtain a systematic computer-assisted method 
for computing the polynomial expressions of $x_{14,15}, \cdots, x_{45}$ as elements of $R$. 

As explained in the proof of \cite[Theorem 4.8]{Nagasato5}, the same polynomial expressions 
of $s_{ak}$ can alternatively be obtained by first sliding the subarc $c_k$ 
all the way to the right-hand side of $D$ in Figure \ref{diagram_T45}, 
and then resolving the winding parts at the end. 
(For details of this procedure, see the proof of Theorem 4.8 in \cite{Nagasato5}.)   

We describe the above diagrammatic elimination process case by case for each parameter 
$x_{ij}$ $(1 \leq i \leq j \leq 15)$:
\begin{enumerate}
\item
For $x_{ij}$ $(1 \leq i \leq j \leq 4)$, place the corresponding loop $s_{ij}$ in $D$ 
so that the subarc $c_i$ (resp.\ $c_j$) hooks the $i$th (resp.\ $j$th) arc 
on the left-hand side of $D$, while the other subarc $c_j$ (resp.\ $c_i$) lies on 
the right-hand side of $D$. Slide the subarc $c_i$ (resp. $c_j$) to the right-hand side of $D$, 
keeping the other subarc $c_j$ (resp.\ $c_i$) and the endpoints of $c_i$ 
(resp.\ $c_j$) fixed. Then resolve any resulting winding parts at the end. 
After substituting $s_{ij}=x_{ij}$, the resulting expression yields a polynomial in $R$, 
which we denote by $g_i(x_{ij})$ (resp.\ $g_j(x_{ij})$). 

\item
For $x_{ij}$ $(1 \leq i < j \leq 4)$ in the additional relations, 
place the corresponding loop $s_{ij}$ on the left-hand side of $D$ 
and slide it to the right-hand side of $D$.
In this case, the resulting loop has no winding parts in $D$.
Hence we interpret the resulting loop directly as an element of $R$, after substituting 
$s_{ij}=x_{ij}$. 

\item
For the remaining $x_{ij}$ $(i<j)$, 
place the corresponding loop $s_{ij}$ in $D$ and slide it to the right-hand side of $D$.
Then resolve the winding parts of the resulting loop at the end. 
After substituting $s_{ij}=x_{ij}$, the resulting expression yields a polynomial in $R$.
\end{enumerate}
Process (1) yields the relations $x_{ij}=g_i(x_{ij})$ and $x_{ij}=g_j(x_{ij})$ 
for all $1 \leq i \leq j \leq 4$. 
Process (2) implies that $x_{12}=x_{23}=x_{34}=x_{14}$ and $x_{13}=x_{24}$. 
Process (3) provides expressions for $x_{ij}$ (with $5 \leq i$ or $5 \leq j$) 
as polynomials in $R$, thereby eliminating these parameters. 
Note that in Processes (2) and (3), the case of $x_{ii}$ $(1 \leq i \leq 15)$ are omitted, 
since they yield only the trivial relations $x_{ii}=2$. 

By the above description, we define a biregular map (that is, an isomorphism onto its image) 
$i: F_2(T_{4,5}) \to \mathrm{Im}(i) \subset \C^{4 \choose 2}$ by 
\[
(x_{12},\cdots,x_{14,15}) \mapsto (x_{12},x_{13},x_{14},x_{23},x_{24},x_{34}). 
\]
The relations $x_{ij}=g_i(x_{ij})$ and $x_{ij}=g_j(x_{ij})$ in Process (1), 
and $x_{12}=x_{23}=x_{34}=x_{45}$ and $x_{13}=x_{24}$ in Process (2)  
become the defining polynomials of $\mathrm{Im}(i)$. Thus, 
\[
F_2(T_{4,5}) \cong \left\{(x_{12},\cdots,x_{34})\in \C^{4 \choose 2}\ \left|\ 
\begin{array}{l}
\mbox{$x_{ij}=g_{i}(x_{ij})$, $x_{ij}=g_{j}(x_{ij})$ $(1 \leq i \leq j \leq 4)$}\\
\mbox{$x_{12}=x_{23}=x_{34}=x_{14}$, $x_{13}=x_{24}$}
\end{array}\right.\right\}, 
\]
whose parametrization corresponds to the presentation in (\ref{t45-4bridge}). 
Moreover, we can eliminate $x_{14}$, $x_{23}$, $x_{24}$, and $x_{34}$ using the relations 
$x_{12}=x_{23}=x_{34}=x_{14}$ and $x_{13}=x_{24}$. 
Then the diagrammatic arguments above implies that, the relations obtained from 
$x_{ij}=g_{i}(x_{ij})$ and $x_{ij}=g_{j}(x_{ij})$ by substituting 
\[
x_{12}=x_{23},\ x_{23}=x_{34},\ x_{34}=x_{14},\ x_{14}=x_{12},\ x_{13}=x_{24},\ x_{24}=x_{13}
\]
are respectively 
\[
x_{i+1,j+1}=g_{i+1}(x_{i+1,j+1}),\ x_{i+1,j+1}=g_{j+1}(x_{i+1,j+1}), 
\]
where the indices are taken cyclically from $1$ to $4$, that is, 
if $i$ (resp. $j$) is $4$, then $i+1=1$ (resp. $j+1=1$). 
Consequently, the relations $x_{ij}=g_{i}(x_{ij})$ and $x_{ij}=g_{j}(x_{ij})$ 
for $(i,j)=(1,4)$, $(2,3)$, and $(3,4)$ reduce to 
\[
x_{12}=g_{1}(x_{12}),\ x_{12}=g_{2}(x_{12}),
\] 
while those for $(i,j)=(2,4)$ reduce to
\[
x_{13}=g_{1}(x_{13}),\ x_{13}=g_{3}(x_{13}). 
\] 
Moreover, the relations $x_{ii}=g_i(x_{ii})$ for $2 \le i \le 4$ reduce to
\[
x_{11}=g_1(x_{11}).
\]
Therefore, $F_2(T_{4,5})$ is isomorphic to the following algebraic set: 
\[
F_2(T_{4,5}) \cong \left\{(x_{12},x_{13}) \in \C^2\ \left|\ 
\begin{array}{l}
x_{1j}=\widetilde{g_{1}}(x_{1j}),\ x_{1j}=\widetilde{g_j}(x_{1j})\ (j=2,3)\\
x_{11}=\widetilde{g_1}(x_{11})
\end{array}\right.\right\}, 
\]
where $\widetilde{g_{i}}(x_{ij})$ (resp. $\widetilde{g_j}(x_{ij})$) denotes the polynomial 
obtained from $g_{i}(x_{ij})$ (resp. $g_j(x_{ij})$) by substituting 
\[
x_{14}=x_{12},\ x_{23}=x_{12},\ x_{24}=x_{13},\ x_{34}=x_{12}. 
\]
We can compute the following explicit expressions of $\widetilde{g_{i}}(x_{ij})$ 
by the above diagrammatic argument\footnote{The polynomials 
$\widetilde{g_{*}}(x_{ij})$ are computed using a computer program written 
to implement the action of the braid group on $R$ described above. 
The same program will be also applied to the case of $T_{5,6}$.}. 
Set $a=x_{12}$ and $b=x_{13}$. 
\begin{eqnarray*}
x_{11}=2=\widetilde{g_{1}}(x_{11}) &=& a^5-4a^3b+3a^3+3ab^2-2ab-3a,\\
a=\widetilde{g_{1}}(x_{12}) &=& a^6-4a^4b+2a^4+3a^2b^2+a^2b-5a^2-b^2+2,\\
a=\widetilde{g_{2}}(x_{12}) &=& a^4b-a^4-3a^2b^2+4a^2b+b^3-3b,\\
b=\widetilde{g_{1}}(x_{13}) &=& a^5b-a^5-4a^3b^2+6a^3b+3ab^3-a^3-3ab^2-5ab+3a,\\
b=\widetilde{g_{3}}(x_{13}) &=& a^5-3a^3b+a^3+ab^2+2ab-3a. 
\end{eqnarray*}
Solving these equations, we see that $F_2(T_{4,5}) \subset \C^2$ 
consists of the following 6 points
\[
(x_{12},x_{13})=(2,2),\ (-1,1),\ (\alpha,-2+2\alpha),\ (\beta,2),
\]
where $\alpha$ and $\beta$ are the roots of the polynomials
$z^2-3z+1$ and $z^2+z-1$, respectively.

We show that the point $(x_{12},x_{13})=(-1,1)$ is a ghost character of $T_{4,5}$. 
As described in the proof of Theorem 4.8 in \cite{Nagasato5}, 
the above elimination process reduce each of the parameters $x_{ijk}$ $(1 \leq i<j<k \leq 15)$ 
for $S_0(T_{4,5}) \subset \C^{{15 \choose 2}+{15 \choose 3}}$ to a linear combinations of 
\[
x_{123},\ x_{124},\ x_{134},\ x_{234}.
\] 
This reduction provides a biregular map 
$j: S_0(T_{4,5}) \to {\rm Im}(j) \subset i(F_2(T_{4,5})) \times \C^4$ 
defined by 
\[
(x_{12},\cdots,x_{14,15};x_{123},\cdots,x_{13,14,15}) 
\mapsto 
(x_{12},\cdots,x_{45};x_{123},x_{124},x_{134},x_{234}). 
\] 
Then the relations (GH) for $S_0(T_{4,5}) \subset \C^{{15  \choose 2}+{15 \choose 3}}$ 
reduce to those for $j(S_0(T_{4,5}))$ in $\C^{10}$. 
Here, the point $(x_{12},x_{13})=(-1,1) \in F_2(T_{4,5}) \subset \C^2$ corresponds to 
\[
\mathbf{x}=(x_{12},x_{13},x_{14},x_{23},x_{24},x_{34})=(-1,1,-1,-1,1,-1)
\in i(F_2(T_{4,5})) \subset \C^{4 \choose 2}.
\]
A direct computation shows that the point $\mathbf{x}$ fails to satisfy 
one of the reduced relations (GH) for $j(S_0(T_{4,5}))$. 
To see this, one may check that $\mathbf{x}$ does not satisfy 
the rectangle relation (R) for $j(S_0(T_{4,5}))$: 
\[
\left|
\begin{array}{cccc}
2 & x_{12} & x_{13} & x_{14}\\
x_{21} & 2 & x_{23} & x_{24}\\
x_{31} & x_{32} & 2 & x_{34}\\
x_{41} & x_{42} & x_{43} & 2
\end{array}
\right|=0.
\]
We note that any point satisfying the relations (GH) 
also satisfies the relations (R) (see the proof of Proposition 3.1 in \cite{Nagasato5}). 
Therefore, $\mathbf{x}$ fails to satisfy at least one of the reduced relations (GH) 
for $j(S_0(T_{4,5}))$, and hence is a ghost character of $T_{4,5}$. 

On the other hand, one can verify that all other points in $F_2(T_{4,5}) \subset \C^2$ 
satisfy the reduced relations (GH) for $j(S_0(T_{4,5}))$. 
Consequently, we obtain the following. 

\begin{proposition}
The $(4,5)$-torus knot $T_{4,5}$ admits exactly one ghost character, 
given by $(x_{12},x_{13})=(-1,1) \in F_2(T_{4,5}) \subset \C^2$. 
\end{proposition}

Combined with Theorem 4.10 in \cite{Nagasato5}, this proposition shows that 
$T_{4,5}$ provides a counterexample to Ng's conjecture. 
Indeed, in the following section, we will observe that $h^*$ is surjective but not injective. 
Moreover, this behavior of $h^*$ shows that $T_{4,5}$ provides the first example 
for which the surjectivity of the map $\widehat{\Phi}$ fails 
(see Theorem 4.9 (1) and (2) in \cite{Nagasato5}). 


\subsection{Ghost characters of $T_{5,6}$}\label{subsec_ghost_t56}
We now turn to the computation of ghost characters of the $(5,6)$-torus knot $T_{5,6}$. 
Since the strategy is the same as that for $T_{4,5}$, 
we omit routine details and focus on the points where the computation differs. 

\begin{figure}[htbp]
\[
D=\begin{minipage}{11cm}
\begin{overpic}[width=\hsize]{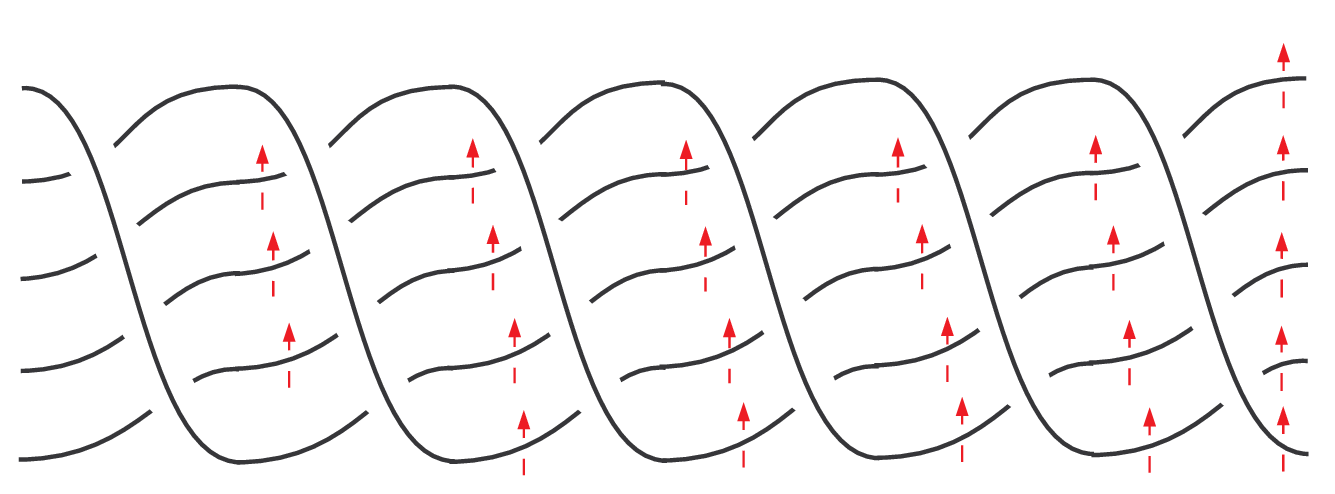}
\put(12,26){$m_{24}$}
\put(13,19){$m_{23}$}
\put(14,12){$m_{22}$}
\put(28,26){$m_{21}$}
\put(29,19){$m_{20}$}
\put(31,12){$m_{19}$}
\put(32,5){$m_{18}$}
\put(44,26){$m_{17}$}
\put(45,19){$m_{16}$}
\put(47,12){$m_{15}$}
\put(48,5){$m_{14}$}
\put(60,26){$m_{13}$}
\put(61,19){$m_{12}$}
\put(63,12){$m_{11}$}
\put(64,5){$m_{10}$}
\put(76,26){$m_9$}
\put(77,19){$m_8$}
\put(79,12){$m_7$}
\put(80,5){$m_6$}
\put(98,32){$m_5$}
\put(98,25){$m_4$}
\put(98,18){$m_3$}
\put(98,11){$m_2$}
\put(98,5){$m_1$}
\end{overpic}
\end{minipage}
\]
\caption{A diagram of $T_{5,6}$ in braid (bridge) position and meridians $m_1,\cdots,m_{24}$. 
The five parallel curves connecting the both sides of the diagram are omitted.}
\label{diagram_T56}
\end{figure}

Let $D$ be the diagram of the $(5,6)$-torus knot $T_{5,6}$ shown in Figure \ref{diagram_T56}. 
We assign meridians $m_1,\cdots,m_{24}$ as indicated in Figure $\ref{diagram_T56}$, 
which yields the Wirtinger presentation of $G(T_{5,6})$ associated with $D$. 

With this notation, the fundamental variety $F_2(T_{5,6})$ is given by 
\[
F_2(T_{5,6})=\left\{(x_{12},\cdots,x_{23,24}) \in \C^{24 \choose 2}\ \left|\ 
\begin{array}{c}
x_{ak}=x_{ij}x_{ai}-x_{aj} \mbox{ for any $1 \leq a \leq 24$}\\
\mbox{and any Wirtinger triple $(i,j,k)$}
\end{array}
\right.\right\}. 
\]
As in the case of $T_{4,5}$, we apply the elimination process based on 
the braid position of $T_{5,6}$ to obtain a reduced description of $F_2(T_{5,6})$.
We begin with the fundamental relations (F2). For $1 \leq a \leq 24$, they are given by 
\[
\begin{array}{ll}
x_{a,24}=x_{14,18}x_{a,18}-x_{a,14}, & x_{a,23}=x_{18,21}x_{a,18}-x_{a,21},\\
x_{a,22}=x_{18,20}x_{a,18}-x_{a,20}, & x_{a,21}=x_{10,14}x_{a,14}-x_{a,10},\\
x_{a,20}=x_{14,17}x_{a,14}-x_{a,17}, & x_{a,19}=x_{14,16}x_{a,14}-x_{a,16},\\
x_{a,18}=x_{14,15}x_{a,14}-x_{a,15}, & x_{a,17}=x_{6,10}x_{a,10}-x_{a,6},\\
x_{a,16}=x_{10,13}x_{a,10}-x_{a,13}, & x_{a,15}=x_{10,12}x_{a,10}-x_{a,12},\\
x_{a,14}=x_{10,11}x_{a,10}-x_{a,11}, & x_{a,13}=x_{1,6}x_{a,6}-x_{a,1},\\
x_{a,12}=x_{6,9}x_{a,6}-x_{a,9}, & x_{a,11}=x_{6,8}x_{a,6}-x_{a,8},\\
x_{a,10}=x_{6,7}x_{a,6}-x_{a,7}, & x_{a,9}=x_{1,5}x_{a,1}-x_{a,5},\\
x_{a,8}=x_{1,4}x_{a,1}-x_{a,4},  & x_{a,7}=x_{1,3}x_{a,1}-x_{a,3},\\
x_{a,6}=x_{1,2}x_{a,1}-x_{a,2}, & \\
x_{a,5}=x_{18,19}x_{a,18}-x_{a,19}, & x_{a,4}=x_{5,18}x_{a,5}-x_{a,18}\\
x_{a,3}=x_{5,24}x_{a,5}-x_{a,24}, & x_{a,2}=x_{5,23}x_{a,5}-x_{a,23},\\
x_{a,1}=x_{5,22}x_{a,5}-x_{a,22}. &
\end{array}
\]
The last 5 types of (F2) are written for Wirtinger triples $(i,j,k)$ 
with $j > k$, using the symmetry between $j$ and $k$, for technical reasons. 
Similar to the case of $T_{4,5}$, we show that the elimination process can be carried out 
mostly manually. 

Analogously, to facilitate the process, we introduce the following additional relations 
derived from (F2): 
\[
\begin{array}{ll}
\fbox{$x_{12}$} &=x_{5,22}x_{25}-x_{2,22}=x_{5,22}x_{5,23}-(x_{5,23}x_{5,22}-x_{22,23})
=\fbox{$x_{22,23}$},\\
\fbox{$x_{13}$} &=x_{5,24}x_{15}-x_{1,24}=x_{5,24}x_{5,22}-(x_{5,22}x_{5,24}-x_{22,24})
=\fbox{$x_{22,24}$},\\
\fbox{$x_{14}$} &=x_{5,18}x_{15}-x_{1,18}=x_{5,18}x_{5,22}-(x_{5,22}x_{5,18}-x_{18,22})
=\fbox{$x_{18,22}$},\\
\fbox{$x_{15}$} &=x_{5,22}=\fbox{$x_{18,19}x_{18,22}-x_{19,22}$},\\
\fbox{$x_{23}$} &=x_{5,23}x_{35}-x_{3,23}=x_{5,23}x_{5,24}-(x_{5,24}x_{5,23}-x_{23,24})
=\fbox{$x_{23,24}$},\\
\fbox{$x_{24}$} &=x_{5,23}x_{45}-x_{4,23}=x_{5,23}x_{5,18}-(x_{5,18}x_{5,23}-x_{18,23})
=\fbox{$x_{18,23}$},\\
\fbox{$x_{25}$} &=x_{5,23}=\fbox{$x_{18,19}x_{18,23}-x_{19,23}$},\\
\fbox{$x_{34}$} &=x_{5,24}x_{45}-x_{4,24}=x_{5,24}x_{5,18}-(x_{5,18}x_{5,24}-x_{18,24})
=\fbox{$x_{18,24}$},\\
\fbox{$x_{35}$} &=x_{5,24}=\fbox{$x_{18,19}x_{18,24}-x_{19,24}$},\\
\fbox{$x_{45}$} &=x_{5,18}=\fbox{$x_{18,19}$}.
\end{array}
\]
The elimination procedure proceeds exactly as in the case of $T_{4,5}$. 
We first eliminate $x_{a,24}$ $(1 \leq a \leq 23)$ from (F2), using the relations 
\[
x_{a,24}=x_{14,18}x_{a,18}-x_{a,14}.
\] 
Next, we eliminate $x_{a,23}$ $(1 \leq a \leq 22)$ from (F2), 
including the equations obtained in the previous step, using 
\[
x_{a,23} = x_{18,21}x_{a,18}-x_{a,21},
\]
and then continue recursively down to $x_{a,6}$. 
After this process, the parameters $x_{a,24}$ ($1 \leq a \leq 23$), $\cdots$, 
$x_{a,6}$ $(1 \leq a \leq 5)$ are expressed as polynomials in the ring 
\[
\C[x_{12},x_{13},x_{14},x_{15},x_{23},x_{24},x_{25},x_{34},x_{35},x_{45}].
\]

By the elimination procedure analogous to the case of $T_{4,5}$, we obtain a biregular map 
(i.e., an isomorphism onto its image) 
$i: F_2(T_{5,6}) \to \mathrm{Im}(i) \subset \C^{5 \choose 2}$ defined by  
\[
(x_{12},\cdots,x_{23,24}) \mapsto (x_{12},x_{13},x_{14},x_{15},x_{23},x_{24},x_{25},x_{34},x_{35},x_{45}).
\] 
The defining equations of $\mathrm{Im}(i)$ consist of the relations 
$x_{ij}=g_i(x_{ij})$ and $x_{ij}=g_j(x_{ij})$ for $1 \le i < j \le 5$ obtained from Process (1), 
together with $x_{12}=x_{23}=x_{34}=x_{45}=x_{15}$ and $x_{13}=x_{24}=x_{35}=x_{14}=x_{25}$ 
obtained from Process (2). Hence we obtain the following: 
\[
F_2(T_{5,6}) \cong \left\{(x_{12},\cdots,x_{45})\in \C^{5 \choose 2}\ \left|\ 
\begin{array}{l}
\mbox{$x_{ij}=g_{i}(x_{ij})$, $x_{ij}=g_{j}(x_{ij})$ $(1 \leq i \leq j \leq 5)$}\\
\mbox{$x_{12}=x_{23}=x_{34}=x_{45}=x_{15}$},\\
\mbox{$x_{13}=x_{24}=x_{35}=x_{14}=x_{25}$}
\end{array}\right.\right\}. 
\]
The resulting parametrization of $F_2(T_{5,6})$ corresponds to the following $5$-bridge knot group 
presentation of $T_{5,6}$ (see the presentation in (\ref{t56-5bridge}) for more precise description): 
\begin{eqnarray*}
G(T_{5,6}) \cong \langle m_1,\ldots,m_5 \mid w_1,\cdots,w_4 \rangle.
\end{eqnarray*}

We now eliminate $x_{14},x_{15},x_{23},x_{24},x_{25},x_{34},x_{35},x_{45}$ using the relations
\[
x_{12}=x_{23}=x_{34}=x_{45}=x_{15},\ x_{13}=x_{24}=x_{35}=x_{14}=x_{25}. 
\]
Moreover, these relations yield 
\begin{eqnarray*}
&&x_{12}=x_{23},\ x_{23}=x_{34},\ x_{34}=x_{45},\ x_{45}=x_{15},\ x_{15}=x_{12},\\
&&x_{13}=x_{24},\ x_{24}=x_{35},\ x_{35}=x_{14},\ x_{14}=x_{25},\ x_{25}=x_{13}. 
\end{eqnarray*}
By substituting these equations into $x_{ij}=g_{i}(x_{ij})$ and $x_{ij}=g_{j}(x_{ij})$, 
we obtain $x_{i+1,j+1}=g_{i+1}(x_{i+1,j+1})$ and $x_{i+1,j+1}=g_{j+1}(x_{i+1,j+1})$, respectively, 
where indices are taken cyclically from $1$ to $5$. 
Hence $x_{ij}=g_{i}(x_{ij})$, $x_{ij}=g_{j}(x_{ij})$ for $(i,j)=(1,5),(2,3),(3,4),(4,5)$ reduce to 
\[
x_{12}=g_{1}(x_{12}),\ x_{12}=g_{2}(x_{12}),
\] 
and those for $(i,j)=(1,4),(2,4),(2,5),(3,5)$ reduce to 
\[
x_{13}=g_{1}(x_{13}),\ x_{13}=g_{3}(x_{13}), 
\]
via this substitution. 
Furthermore, $x_{ii}=g_{i}(x_{ii})$ for $2 \leq i \leq 5$ reduce to
\[
 x_{11}=g_{1}(x_{11}). 
\]
Therefore, $F_2(T_{5,6})$ is isomorphic to the following algebraic set:  
\[
F_2(T_{5,6}) \cong \left\{(x_{12},x_{13}) \in \C^2\ \left|\ 
\begin{array}{l}
x_{1j}=\widetilde{g_1}(x_{1j}),\ x_{1j}=\widetilde{g_j}(x_{1j})\ (j=2,3)\\
x_{11}=\widetilde{g_1}(x_{11})
\end{array}\right.\right\}, 
\]
where $\widetilde{g_{i}}(x_{ij})$ (resp. $\widetilde{g_j}(x_{ij})$) denotes the polynomial 
obtained from $g_{i}(x_{ij})$ (resp. $g_j(x_{ij})$) by substituting 
$x_{14}=x_{13}$, $x_{15}=x_{12}$, $x_{23}=x_{12}$, $x_{24}=x_{13}$, $x_{25}=x_{13}$,  
$x_{34}=x_{12}$, $x_{35}=x_{13}$, $x_{45}=x_{12}$. 
Using the computer-assisted method developed in Subsection \ref{subsec_ghost_t45}, 
we compute the explicit expressions for $\widetilde{g_i}(x_{ij})$. 
Let $a=x_{12}$, $b=x_{13}$. 
\begin{eqnarray*}
2=\widetilde{g_{1}}(x_{11})&=& a^6-5a^4b+4a^3b+6a^2b^2-3a^3-6ab^2-b^3+2ab+b^2+3a,\\
a=\widetilde{g_{1}}(x_{12})&=& a^7-5a^5b-a^5+4a^4b+6a^3b^2-3a^4+4a^3b-6a^2b^2-ab^3-a^2b\\
&&-2ab^2+5a^2+2b^2-2,\\
a=\widetilde{g_{2}}(x_{12})&=& a^5b-a^4b-4a^3b^2+a^4+6a^2b^2+3ab^3-4a^2b-2ab^2-3b^3-a^2\\
&&+2ab+3b,\\
b=\widetilde{g_{1}}(x_{13})
&=& a^6b-a^6-5a^4b^2+4a^4b+4a^3b^2+6a^2b^3+a^4-6a^3b-3a^2b^2-6ab^3\\
&&-b^4+2a^3-3a^2b+4ab^2+b^3+5ab+b^2-3a,\\
b=\widetilde{g_{3}}(x_{13})&=&a^5b-a^5-4a^3b^2+3a^3b+3a^2b^2+3ab^3+a^3-4a^2b-ab^2-2b^3\\
&&+a^2-2ab+3b. 
\end{eqnarray*}
By solving these equations, we see that $F_2(T_{5,6}) \subset \C^2$ consists of 
the following 10 points 
\[
(x_{12},x_{13})=(2,2),\ (0,-1),\ (1,1),\ (-2,1),\ (z_1,-1+2z_1),\ (z_2,1),\ (z_3,-1-z_3),
\]
where $z_1$, $z_2$, and $z_3$ are the roots of the equations
$z^2-5z+5=0$, $z^2-z-1=0$, and $z^2+z-1=0$, respectively.

We show that the points $(x_{12},x_{13})=(0,-1),(1,1),(-2,1)$ in $F_2(T_{5,6})$ 
are the ghost characters of $T_{5,6}$. 
Analogously to the case of $T_{4,5}$, each of the parameters $x_{ijk}$ $(1 \leq i<j<k \leq 24)$ 
for $S_0(T_{5,6}) \subset \C^{{24 \choose 2}+{24 \choose 3}}$ 
reduces to a linear combination of 
\[
x_{123},\ x_{124},\ x_{125},\ x_{134},\ x_{135},\ x_{145},\ x_{234},\ x_{235},\ x_{245},\ x_{345}. 
\]
This reduction provides a biregular map 
$j: S_0(T_{5,6}) \to \mathrm{Im}(j) \subset i(F_2(T_{5,6})) \times \C^{10}$ 
defined by 
\[
(x_{12},\cdots,x_{23,24};x_{123},\cdots,x_{22,23,24}) 
\mapsto 
(x_{12},\cdots,x_{45};x_{123},\cdots,x_{345}). 
\]
Then the relations (GH) for $S_0(T_{5,6}) \subset \C^{{5 \choose 2}+{5 \choose 3}}$ 
reduce to those for $j(S_0(T_{5,6})) \subset \C^{20}$.  
As demonstrated for the case $T_{4,5}$, the points in $i(F_2(T_{5,6}))$ corresponding to 
$(x_{12},x_{13})=(0,-1),(1,1),(-2,1)$ fail to satisfy 
one of the reduced relations (GH) for $j(S_0(T_{5,6}))$. 
To see this, one may check that the points do not satisfy 
the following relation (R) for $j(S_0(T_{5,6}))$: 
\[
\left|
\begin{array}{cccc}
2 & x_{12} & x_{13} & x_{14}\\
x_{21} & 2 & x_{23} & x_{24}\\
x_{31} & x_{32} & 2 & x_{34}\\
x_{41} & x_{42} & x_{43} & 2
\end{array}
\right|=0. 
\]
Therefore, the points $(x_{12},x_{13})=(0,-1), (1,1), (-2,1) \in F_2(T_{5,6})$ are 
ghost characters of $T_{5,6}$. 
A direct calculation shows that the remaining points in $F_2(T_{5,6})$ 
satisfy all reduced relations (GH) for $j(S_0(T_{5,6}))$. Consequently, we obtain the following. 

\begin{proposition}\label{prop_ghost}
The $(5,6)$-torus knot $T_{5,6}$ admits exactly three ghost characters, given by 
$(x_{12},x_{13})=(0,-1)$, $(1,1)$, $(-2,1) \in F_2(T_{5,6}) \subset \C^2$.
\end{proposition}

Combined with \cite[Theorem 4.10]{Nagasato5}, this proposition proves that 
$T_{5,6}$ provides a counterexample to Ng's conjecture. 
Indeed, in the next section, we will observe that the map $h^*$ is neither surjective nor injective. 


\section{Behavior of the map $h^*$}\label{sec_h}
As observed above, the map $h^* : X(\Sigma_2K) \to F_2(K)$ is not an isomorphism 
for torus knots $T_{4,5}$ and $T_{5,6}$. 
For these knots, we investigate the behavior of the map $h^*$ 
with respect to surjectivity and injectivity, using ghost characters. 
The computations in this section are based on the master's thesis of the second author \cite{Suzuki}. 
Following the strategy outlined at the end of Section \ref{overview}, we identify ghost characters 
that obstruct either the surjectivity or the injectivity of $h^*$. 


\subsection{Character variety $X(\Sigma_2T_{4,5})$ and behavior of $h^*$}
As shown in Subsection \ref{subsec_ghost_t45}, the $(4,5)$-torus knot $T_{4,5}$ admits 
a unique ghost character 
\[
\mathbf{g}:=(x_{12},x_{13})=(-1,1) \in F_2(T_{4,5}) \subset \C^2.
\] 
Hence, 
if there exists a representation $\rho_*: \pi_1(\Sigma_2T_{4,5}) \to \SL_2(\C)$ satisfying 
\[
t_{m_1m_2}(\chi_{\rho_*})=-1,\ t_{m_1m_3}(\chi_{\rho_*})=1,
\]
then the map $h^*$ is surjective but not injective, as shown in the proof of Theorem 4.10 
in \cite{Nagasato5}, whereas the map $\widehat{\Phi}$ fails to be surjective 
by \cite[Theorem 4.9 (2)]{Nagasato5}. 
Indeed, such a representation exists, leading to the following result.  

\begin{theorem}\label{thm_main2}
For the $(4,5)$-torus knot $T_{4,5}$, the map $h^*$ is surjective but not injective, 
whereas $\widehat{\Phi}$ fails to be surjective. 
In other words, there exists an $\SL_2(\C)$-representation of $\pi_1(\Sigma_2T_{4,5})$ 
that cannot be obtained, via the map $\widehat{\Phi}$ (in the sense of $(\ref{rep_s2k})$), 
from any trace-free $\SL_2(\C)$-representation of the knot group $G(T_{4,5})$. 
\end{theorem}

The knot $T_{4,5}$ provides the first example for which the surjectivity of 
$\widehat{\Phi}$ fails. 

\begin{proof}
We begin by computing $\pi_1(\Sigma_2T_{4,5})$. 

\begin{lemma}\label{lem_pi_1s2k}
Let $G(T_{4,5})=\langle m_1,\cdots m_{15} \mid r_1,\cdots, r_{14} \rangle$ 
be the Wirtinger presentation\footnote{The last relator $r_{15}$ is eliminated 
due to deficiency one property of knot groups.} 
associated with the diagram $D$ in Figure $\ref{diagram_T45}$.  
Set $x=m_1m_2$, $y=m_1m_3$, and $z=m_1m_4$. 
Then the fundamental group $\pi_1(\Sigma_2T_{4,5})$ admits the presentation  
\[
\pi_1(\Sigma_2T_{4,5}) \cong \langle x,y,z \mid w_i\ (1 \leq i \leq 3),
m_j^2\ (1 \leq j \leq 4 )\rangle, 
\]
where the relators $w_i$ are given by 
\begin{eqnarray*}
w_1 &=& z^{-1} x^{-1} y z^{-1} x z^{-1} y x^{-1} z^{-1},\\
w_2 &=& z^{-1} x^{-1} y z^{-1} y z^{-1} y x^{-1} z^{-1} x,\\
w_3 &=& z^{-1} x^{-1} y z^{-1} y x^{-1} z^{-1} y.
\end{eqnarray*}
\end{lemma}

We explain the computation of $\pi_1(\Sigma_2 T_{4,5})$ in Lemma \ref{lem_pi_1s2k}.  
First, using the relators $r_1,\ldots,r_{14}$ of the Wirtinger presentation of $G(T_{4,5})$, 
we obtain 
\begin{eqnarray*}
m_5    &=& m_1 m_2 m_1^{-1},\\
m_6    &=& m_1 m_3 m_1^{-1},\\
m_7    &=& m_1 m_4 m_1^{-1},\\
m_8    &=& m_5 m_6 m_5^{-1} = m_1 m_2 m_3 m_2^{-1} m_1^{-1},\\
m_9    &=& m_5 m_7 m_5^{-1} = m_1 m_2 m_4 m_2^{-1} m_1^{-1},\\
m_{10} &=& m_5 m_1 m_5^{-1} = m_1 m_2 m_1 m_2^{-1} m_1^{-1},\\
m_{11} &=& m_8 m_9 m_8^{-1} = m_1 m_2 m_3 m_4 m_3^{-1} m_2^{-1} m_1^{-1},\\
m_{12} &=& m_8 m_{10} m_8^{-1} = m_1 m_2 m_3 m_1 m_3^{-1} m_2^{-1} m_1^{-1},\\
m_{13} &=& m_8 m_5 m_8^{-1} = m_1 m_2 m_3 m_2 m_3^{-1} m_2^{-1} m_1^{-1},\\
m_{14} &=& m_{11} m_{13} m_{11}^{-1} = m_1 m_2 m_3 m_4 m_2 m_4^{-1} m_3^{-1} m_2^{-1} m_1^{-1},\\
m_{15} &=& m_{11} m_8 m_{11}^{-1} = m_1 m_2 m_3 m_4 m_3 m_4^{-1} m_3^{-1} m_2^{-1} m_1^{-1},\\
m_1    &=& m_4 m_{14} m_4^{-1} = m_4 m_1 m_2 m_3 m_4 m_2 m_4^{-1} m_3^{-1} m_2^{-1} 
m_1^{-1} m_4^{-1},\\
m_2    &=& m_4 m_{15} m_4^{-1} = m_4 m_1 m_2 m_3 m_4 m_3 m_4^{-1} m_3^{-1} m_2^{-1} 
m_1^{-1} m_4^{-1},\\
m_3    &=& m_4 m_{11} m_4^{-1} = m_4 m_1 m_2 m_3 m_4 m_3^{-1} m_2^{-1} m_1^{-1} m_4^{-1}.
\end{eqnarray*}
By applying Tietze transformations, the first 11 relations show that 
$G(T_{4,5})$ is generated by $m_1$, $m_2$, $m_3$, and $m_4$, 
and that the set of relators of $G(T_{4,5})$ is normally generated by the last 3 relations. 
Hence, we obtain the following 4-bridge knot group presentation: 
\begin{eqnarray}\label{t45-4bridge-2}
G(T_{4,5})=\langle m_1,m_2,m_3,m_4 \mid w_1,w_2,w_3 \rangle, 
\end{eqnarray}
where $w_1$, $w_2$ and $w_3$ are given by 
\begin{eqnarray*}
w_1 &=& m_4 m_1 m_2 m_3 m_4 m_2 m_4^{-1} m_3^{-1} m_2^{-1} m_1^{-1} m_4^{-1} m_1^{-1},\\
w_2 &=& m_4 m_1 m_2 m_3 m_4 m_3 m_4^{-1} m_3^{-1} m_2^{-1} m_1^{-1} m_4^{-1} m_2^{-1},\\
w_3 &=& m_4 m_1 m_2 m_3 m_4 m_3^{-1} m_2^{-1} m_1^{-1} m_4^{-1} m_3^{-1}.
\end{eqnarray*} 
Consequently, by the presentation given in Section \ref{overview}, we obtain  
\begin{eqnarray*}
\pi_1(\Sigma_2 T_{4,5})
& \cong &  \langle m_1m_2, m_1m_3, m_1m_4 
\mid w_i (1 \leq i \leq 6), m_j^2\ (1 \leq j \leq 4 )\rangle, 
\end{eqnarray*}
where $w_4 = m_1 w_1 m_1^{-1}, w_5 = m_1 w_2 m_1^{-1}, w_6 = m_1 w_3 m_1^{-1}$. 
For simplicity, set $x=m_1m_2$, $y=m_1m_3$, and $z=m_1m_4$. Then the relators become 
\begin{eqnarray*}
w_1 &=& z^{-1} x^{-1} y z^{-1} x z^{-1} y x^{-1} z^{-1},\\
w_2 &=& z^{-1} x^{-1} y z^{-1} y z^{-1} y x^{-1} z^{-1} x,\\
w_3 &=& z^{-1} x^{-1} y z^{-1} y x^{-1} z^{-1} y,\\
w_4 &=& z x y^{-1} z x^{-1} z y^{-1} x z,\\
w_5 &=& z x y^{-1} z y^{-1} z y^{-1} x z x^{-1},\\
w_6 &=& z x y^{-1} z y^{-1} x z y^{-1}. 
\end{eqnarray*}
We remark that $w_4 = w_1^{-1}$, $w_5 = x w_2^{-1} x^{-1}$, and $w_6 = y w_3^{-1} y^{-1}$. 
Hence, by applying Tietze transformations, the relators $w_4$, $w_5$, and $w_6$ 
may be omitted. This shows Lemma \ref{lem_pi_1s2k}. 

To complete the proof of Theorem \ref{thm_main2}, we construct, for the ghost character
\[
\mathbf{g}=(x_{12},x_{13})=(-1,1) \in F_2(T_{4,5}) \subset \C^2,
\]
an $\SL_2(\C)$-representation $\rho_*: \pi_1(\Sigma_2T_{4,5}) \to \SL_2(\C)$ 
satisfying
\[
\tr(\rho_*(m_1m_2))=\tr(\rho_*(m_1m_4))=-1,\ \tr(\rho_*(m_1m_3))=1.
\]
Solving the relations in $\pi_1(\Sigma_2T_{4,5})$, we obtain two such representations 
\[
\rho_*=\rho_{\alpha_i}: \pi_1(\Sigma_2 T_{4,5}) \to \SL_2(\C) 
\] 
for $i=1,2$ defined by 
\begin{eqnarray*}
&&(\rho_{\alpha_i}(m_1m_2),\rho_{\alpha_i}(m_1m_3),\rho_{\alpha_i}(m_1m_4))\\
&&=\left(
\left(\begin{array}{cc}e^{\frac{2}{3}\pi \i} & 0\\ 0 & e^{-\frac{2}{3}\pi \i} \end{array}\right),
\left(\begin{array}{cc}-\frac{\i}{\sqrt{3}}e^{\frac{\pi}{3}\i} & -\frac{2}{3}\\
 1 & \frac{\i}{\sqrt{3}}e^{-\frac{\pi}{3}\i} \end{array}\right),
\left(\begin{array}{cc}\frac{\i}{\sqrt{3}}e^{\frac{\pi}{3}\i} & \frac{1+2\alpha_i}{3}\\
\alpha_i & -\frac{\i}{\sqrt{3}}e^{-\frac{\pi}{3}\i} \end{array}\right) 
\right), 
\end{eqnarray*}
where $\alpha_1$ and $\alpha_2$ are the roots of the equation $2x^2+x+2=0$ and $\i:=\sqrt{-1}$.  
Note that $\rho_{\alpha_i}$ is irreducible for $i=1,2$ (see the proof of \cite[Theorem 4.10]{Nagasato5}). 
A direct computation shows   
\begin{eqnarray*}
&& \tr(\rho_{\alpha_i}(m_1m_2))=2\cos\left(\frac{2}{3}\pi\right)=-1,\\
&& \tr(\rho_{\alpha_i}(m_1m_3))=-\frac{\i}{\sqrt{3}} \cdot 2\i \sin\left(\frac{\pi}{3}\right)=1,\\
&&\tr(\rho_{\alpha_i}(m_1m_4))=\frac{\i}{\sqrt{3}} \cdot 2\i \sin\left(\frac{\pi}{3}\right)=-1,\\
&&\tr(\rho_{\alpha_i}(m_1m_2\cdot m_1m_3))
=\tr(\rho_{\alpha_i}(m_1m_2\cdot m_1m_4))
=\tr(\rho_{\alpha_i}(m_1m_3 \cdot m_1m_4))=0,\\
&&\tr(\rho_{\alpha_i}(m_1m_2 \cdot m_1m_3 \cdot m_1m_4))=-\frac{1}{2}-\frac{\sqrt{3}(1+4\alpha_i)}{6}\i.
\end{eqnarray*}
The value $\tr(\rho_{\alpha_i}(m_1m_2 \cdot m_1m_3 \cdot m_1m_4))$ verifies that 
the characters $\chi_{\rho_{\alpha_i}}$ $(i=1,2)$ determine two distinct points in the preimage  
$(h^*)^{-1}(\mathbf{g})$. Thus the map $h^*$ is not injective. 
Moreover, by Theorem 4.9 (2), this implies that $\rho_{\alpha_i}$ $(i=1,2)$ 
cannot be obtained, via $\widehat{\Phi}$ (in the sense of $(\ref{rep_s2k})$), 
from any trace-free $\SL_2(\C)$-representation of $G(T_{4,5})$. 

Regarding the surjectivity of $h^*$, we obtain the following $\SL_2(\C)$-representations of 
$\pi_1(\Sigma_2T_{4,5})$ corresponding to the remaining points of $F_2(T_{4,5})$ through the map $h^*$: 
\[
\begin{array}{lll}
\rho_k(m_1m_2):=
\left(\begin{array}{cc}e^{\frac{2\pi k}{5}\i} & 0\\ 0 & e^{-\frac{2\pi k}{5}\i} \end{array}\right)
& \Rightarrow & 
\tr(\rho_k(m_1m_2))=2\cos\left(\frac{2\pi k}{5}\right)=2\mbox{ or }\frac{-1 \pm \sqrt{5}}{2},\\
\rho_k(m_1m_3):=
\left(\begin{array}{cc}1 & 0\\ 0 & 1 \end{array}\right)
& \Rightarrow & 
\tr(\rho_k(m_1m_3))=2,\\
\rho_k(m_1m_4):=
\left(\begin{array}{cc}e^{\frac{2\pi k}{5}\i} & 0\\ 0 & e^{-\frac{2\pi k}{5}\i} \end{array}\right)
& \Rightarrow &
\tr(\rho_k(m_1m_4))=2\cos\left(\frac{2\pi k}{5}\right)=2\mbox{ or }\frac{-1 \pm \sqrt{5}}{2},\\
\end{array}
\]
where $k=0,1,2,3,4$, and 
\[
\begin{array}{lll}
\rho_{\beta}(m_1m_2):=
\left(\begin{array}{cc} \frac{3 \pm \sqrt{5}+\beta}{4} & 0\\
0 & \frac{3 \pm \sqrt{5}-\beta}{4} \end{array}\right)
& \Rightarrow & 
\tr(\rho_{\beta}(m_1m_2))=\frac{3 \pm \sqrt{5}}{2},\\
\rho_{\beta}(m_1m_3):=
\left(\begin{array}{cc}\frac{66(1 \pm \sqrt{5})+26\beta+\beta^3}{132} & \frac{1 \pm 3\sqrt{5}}{11}\\
1 & \frac{66(1 \pm \sqrt{5})-26\beta-\beta^3}{132} \end{array}\right)
& \Rightarrow & 
\tr(\rho_{\beta}(m_1m_3))=1 \pm \sqrt{5},\\
\rho_{\beta}(m_1m_4):=
\left(\begin{array}{cc}\frac{33(3 \pm \sqrt{5})-7\beta+\beta^3}{132} & \frac{1 \pm 3\sqrt{5}}{11}\\
1 & \frac{33(3 \pm \sqrt{5})+7\beta-\beta^3}{132} \end{array}\right)
& \Rightarrow &
\tr(\rho_{\beta}(m_1m_4))=\frac{3 \pm \sqrt{5}}{2}, 
\end{array}
\]
where $\beta=\sqrt{-2\pm 6\sqrt{5}}$ and the signs in the above descriptions are taken consistently. 
By \cite[Proposition 10]{Nagasato-Yamaguchi}, representations $\rho_k$ and $\rho_{\beta}$ arise from 
irreducible metabelian and irreducible non-metebelian representations of $G(T_{4,5})$, respectively. 
The characters of these representations are mapped under the map $h^*$ 
to all remaining points of $F_2(T_{4,5})$, thereby verifying the surjectivity of the map $h^*$.
This completes the proof. 
\end{proof}

As shown in Theorem 4.8 of \cite{Nagasato5}, all $2$-bridge and $3$-bridge knots 
admit no ghost characters. Hence, by Theorem 4.10 of \cite{Nagasato5}, 
the map $h^*$ is an isomorphism for such a knot. 
Since any torus knot simpler than $T_{4,5}$ has bridge index less than $4$, 
the $(4,5)$-torus knot $T_{4,5}$ is the simplest torus knot for which $h^*$ is not an isomorphism. 
Analogously, by Theorem 4.9 (1) of \cite{Nagasato5}, 
the knot $T_{4,5}$ is the simplest knot for which $\widehat{\Phi}$ is not surjective.


\subsection{Character variety $X(\Sigma_2T_{5,6})$ and behavior of $h^*$}
As shown in Proposition \ref{prop_ghost}, the $(5,6)$-torus knot $T_{5,6}$ 
admits the following 3 ghost characters:
\[
(x_{12},x_{13})=(0,-1),(1,1),(-2,1) \in F_2(T_{5,6}) \subset \C^2.
\]
Accordingly, these 3 points are the candidates for obstructions 
to the surjectivity or injectivity of the map $h^*$. We prove the following result. 

\begin{theorem}\label{thm_main}
The ghost character $\g_1:=(1,1)$ has the property that $(h^*)^{-1}(\g_1)=\emptyset$, 
whereas the ghost character $\g_2:=(0,-1)$ has a preimage $(h^*)^{-1}(\g_2)$ 
consisting of two distinct points. As a consequence, the map $h^*$ is neither 
surjective nor injective.
\end{theorem}

\begin{proof}
We first compute $\pi_1(\Sigma_2T_{5,6})$ using the presentation described
in Section \ref{overview}.

\begin{lemma}[cf.\ Theorem 4.5 in \cite{Nagasato5}]\label{lem_pi_1s2k} 
Let 
\[
G(T_{5,6})=\langle m_1,\ldots,m_{24} \mid r_1,\ldots,r_{23} \rangle
\]
be the Wirtinger presentation associated with the diagram $D$ in Figure $\ref{diagram_T56}$.  
Set $x=m_1m_2$, $y=m_1m_3$, $z=m_1m_4$, and $w=m_1m_5$. 
Then the fundamental group $\pi_1(\Sigma_2T_{5,6})$ admits the presentation 
\[
\pi_1(\Sigma_2T_{5,6}) \cong \langle x,y,z,w \mid w_i\ (1 \leq i \leq 4), 
m_j^2\ (1 \leq j \leq 5)\rangle, 
\]
where the relators $w_i$ $(1 \leq i \leq 4)$ are given by 
\begin{eqnarray*}
w_1&=&w^{-1} x^{-1} y z^{-1} w x^{-1} w z^{-1} y x^{-1} w^{-1},\\
w_2&=&w^{-1} x^{-1} y z^{-1} w y^{-1} w z^{-1} y x^{-1} w^{-1} x,\\
w_3&=&w^{-1} x^{-1} y z^{-1} w z^{-1} w z^{-1} y x^{-1} w^{-1} y,\\
w_4&=&w^{-1} x^{-1} y z^{-1} w z^{-1} y x^{-1} w^{-1} z.
\end{eqnarray*}
\end{lemma}
We now explain how the presentation of $\pi_1(\Sigma_2T_{5,6})$ in Lemma \ref{lem_pi_1s2k} 
is obtained. Using the relators $r_1,\cdots,r_{23}$ of the Wirtinger presentation of $G(T_{5,6})$, 
we obtain the following expressions:  
\begin{eqnarray*}
m_6 &=& m_1 m_2 m_1^{-1},\\
m_7 &=& m_1 m_3 m_1^{-1},\\
m_8 &=& m_1 m_4 m_1^{-1},\\
m_9 &=& m_1 m_5 m_1^{-1},\\
m_{10}&=& m_6 m_7 m_6^{-1}=m_1 m_2 m_3 m_2^{-1} m_1^{-1},\\
m_{11}&=& m_6 m_8 m_6^{-1}=m_1 m_2 m_4 m_2^{-1} m_1^{-1},\\
m_{12}&=& m_6 m_9 m_6^{-1}=m_1 m_2 m_5 m_2^{-1} m_1^{-1},\\
m_{13}&=& m_6 m_1 m_6^{-1}=m_1 m_2 m_1 m_2^{-1} m_1^{-1},\\
m_{14}&=& m_{10}m_{11}m_{10}^{-1}=m_1 m_2 m_3 m_4 m_3^{-1} m_2^{-1} m_1^{-1},\\
m_{15}&=& m_{10}m_{12}m_{10}^{-1}=m_1 m_2 m_3 m_5 m_3^{-1} m_2^{-1} m_1^{-1},\\
m_{16}&=& m_{10}m_{13}m_{10}^{-1}=m_1 m_2 m_3 m_1 m_3^{-1} m_2^{-1} m_1^{-1},\\
m_{17}&=& m_{10}m_6m_{10}^{-1}=m_1 m_2 m_3 m_2 m_3^{-1} m_2^{-1} m_1^{-1},\\
m_{18}&=& m_{14}m_{15}m_{14}^{-1}=m_1m_2m_3m_4m_5m_4^{-1}m_3^{-1}m_2^{-1}m_1^{-1},\\
m_{19}&=& m_{14}m_{16}m_{14}^{-1}=m_1m_2m_3m_4m_1m_4^{-1}m_3^{-1}m_2^{-1}m_1^{-1},\\
m_{20}&=& m_{14}m_{17}m_{14}^{-1}=m_1m_2m_3m_4m_2m_4^{-1}m_3^{-1}m_2^{-1}m_1^{-1},\\
m_{21}&=& m_{14}m_{10}m_{14}^{-1}=m_1m_2m_3m_4m_3m_4^{-1}m_3^{-1}m_2^{-1}m_1^{-1},\\
m_{22}&=& m_{18}m_{20}m_{18}^{-1}=m_1m_2m_3m_4m_5m_2m_5^{-1}m_4^{-1}m_3^{-1}m_2^{-1}m_1^{-1},\\
m_{23}&=& m_{18}m_{21}m_{18}^{-1}=m_1m_2m_3m_4m_5m_3m_5^{-1}m_4^{-1}m_3^{-1}m_2^{-1}m_1^{-1},\\
m_{24}&=& m_{18}m_{14}m_{18}^{-1}=m_1m_2m_3m_4m_5m_4m_5^{-1}m_4^{-1}m_3^{-1}m_2^{-1}m_1^{-1},\\
m_1 &=& m_{5}m_{22}m_{5}^{-1}=m_5m_1m_2m_3m_4m_5m_2m_5^{-1}m_4^{-1}m_3^{-1}m_2^{-1}
m_1^{-1}m_5^{-1},\\
m_2 &=& m_{5}m_{23}m_{5}^{-1}=m_5m_1m_2m_3m_4m_5m_3m_5^{-1}m_4^{-1}m_3^{-1}m_2^{-1}
m_1^{-1}m_5^{-1},\\
m_3 &=& m_{5}m_{24}m_{5}^{-1}=m_5m_1m_2m_3m_4m_5m_4m_5^{-1}m_4^{-1}m_3^{-1}m_2^{-1}
m_1^{-1}m_5^{-1},\\
m_4 &=& m_{5}m_{18}m_{5}^{-1}=m_5m_1m_2m_3m_4m_5m_4^{-1}m_3^{-1}m_2^{-1}m_1^{-1}m_5^{-1}.
\end{eqnarray*}
By applying Tietze transformations, 
the first 19 relations imply that $G(T_{5,6})$ is generated by $m_1$, $m_2$, $m_3$, $m_4$, and $m_5$. 
Moreover, the set of relators of $G(T_{5,6})$ is normally generated by the remaining 4 relations. 
Hence, we obtain the following 5-bridge knot group presentation: 
\begin{eqnarray}\label{t56-5bridge}
G(T_{5,6}) \cong \langle m_1,m_2,m_3,m_4,m_5 \mid w_1,w_2,w_3,w_4 \rangle, 
\end{eqnarray}
where the relators $w_1$, $w_2$, $w_3$, and $w_4$ are given by 
\begin{eqnarray*}
w_1&=&m_5m_1m_2m_3m_4m_5m_2m_5^{-1}m_4^{-1}m_3^{-1}m_2^{-1}m_1^{-1}m_5^{-1}m_1^{-1},\\
w_2&=&m_5m_1m_2m_3m_4m_5m_3m_5^{-1}m_4^{-1}m_3^{-1}m_2^{-1}m_1^{-1}m_5^{-1}m_2^{-1},\\
w_3&=&m_5m_1m_2m_3m_4m_5m_4m_5^{-1}m_4^{-1}m_3^{-1}m_2^{-1}m_1^{-1}m_5^{-1}m_3^{-1},\\
w_4&=&m_5m_1m_2m_3m_4m_5m_4^{-1}m_3^{-1}m_2^{-1}m_1^{-1}m_5^{-1}m_4^{-1}.
\end{eqnarray*}
Then it follows from the presentation given in Section \ref{overview} that 
\begin{eqnarray*}
\pi_1(\Sigma_2 T_{5,6})& \cong &  \langle m_1m_2, m_1m_3, m_1m_4, m_1m_5 \mid w(w_i)\ 
(1 \leq i \leq 4), m_j^2\ (1 \leq j \leq 5) \rangle. 
\end{eqnarray*}
Indeed, starting from the above presentation of $G(T_{5,6})$, we first obtain
\[
\pi_1(\Sigma_2T_{5,6}) \cong \langle m_1m_i\ (2 \leq i \leq 5) \mid 
w(w_j)\ (1 \leq j \leq 8), m_i^2\ (1 \leq i \leq 5) \rangle, 
\] 
where $w_j:=w(m_1w_{j-4}m_1^{-1})$ for $5 \leq j \leq 8$. 
For simplicity, set $x=m_1m_2$, $y=m_1m_3$, $z=m_1m_4$, and $w=m_1m_5$. 
Then, the relators $w_i$ are given by 
\begin{eqnarray*}
w_1&=&w^{-1} x^{-1} y z^{-1} w x^{-1} w z^{-1} y x^{-1} w^{-1},\\
w_2&=&w^{-1} x^{-1} y z^{-1} w y^{-1} w z^{-1} y x^{-1} w^{-1} x,\\
w_3&=&w^{-1} x^{-1} y z^{-1} w z^{-1} w z^{-1} y x^{-1} w^{-1} y,\\
w_4&=&w^{-1} x^{-1} y z^{-1} w z^{-1} y x^{-1} w^{-1} z,\\
w_5&=&w x y^{-1} z w^{-1} x w^{-1} z y^{-1} x w,\\
w_6&=&w x y^{-1} z w^{-1} y w^{-1} z y^{-1} x w x^{-1},\\
w_7&=&w x y^{-1} z w^{-1} z w^{-1} z y^{-1} x w y^{-1},\\
w_8&=&w x y^{-1} z w^{-1} z y^{-1} x w z^{-1}.
\end{eqnarray*}
The relators $w_5,\cdots,w_8$ are conjugate to $w(w_1)^{-1},\cdots,w(w_4)^{-1}$, respectively. 
Therefore, we may omit the relators $w_5,\cdots,w_8$, 
proving Lemma \ref{lem_pi_1s2k}.

To complete the proof of Theorem \ref{thm_main}, we first show that the preimage 
$(h^*)^{-1}(\g_1)$ is empty. Note that $\g_1$ implies $x_{ij}=1$ for all $1 \leq i < j \leq 5$. 
Suppose to the contradiction that there exists an $\SL_2(\C)$-representation $\rho_*$ 
of $\pi_1(\Sigma_2T_{5,6})$ satisfying $t_{m_1m_2}(\rho_*)=1$. 
Since $t_{m_1m_2}(\rho_*) \neq 2$, we may assume, up to conjugation, that 
\begin{eqnarray*}
&&(\rho_*(m_1m_2),\rho_*(m_1m_3),\rho_*(m_1m_4),\rho_*(m_1m_5))\\
&&=\left(
\left(\begin{array}{cc}a & 0\\ 0 & a^{-1} \end{array}\right),
\left(\begin{array}{cc}b & c\\ d & e \end{array}\right),
\left(\begin{array}{cc}f & g\\ h & i \end{array}\right) 
\left(\begin{array}{cc}j & k\\ l & m \end{array}\right) 
\right) \in \SL_2(\C)^4.
\end{eqnarray*}
Computer calculations show that there exists no $\SL_2(\C)$-representation 
$\rho_*$ satisfying
\[
(t_{m_1m_2}(\rho_*),t_{m_1m_3}(\rho_*),t_{m_1m_4}(\rho_*),t_{m_1m_5}(\rho_*))=(1,1,1,1).
\]
To justify this claim, we consider the conditions $t_{m_im_j}(\rho_*)=1$ 
for all $1 \leq i < j \leq 5$. For example, $t_{m_1m_2}(\rho_*)=1$ implies that
\[
a=\frac{1 \pm \sqrt{3}\i}{2}. 
\]
Similarly, the conditions $t_{m_1m_3}(\rho_*)=t_{m_1m_4}(\rho_*)=t_{m_1m_5}(\rho_*)=1$ yield 
\[
e=1-b,\ i=1-f,\ m=1-j.
\]
Since the relation $m_2m_3=(m_1m_2)^{-1}(m_1m_3)$ holds in $\pi_1(\Sigma_2 T_{5,6})$, 
we compute 
\[ 
\tr(\rho_*(m_2m_3)) =\tr\left(\left(\begin{array}{ll}a^{-1} & 0 \\ 0 & a \end{array}\right)
\left(\begin{array}{cc}b & c \\ d & 1-b \end{array}\right)\right) 
=\mp\sqrt{3}\i b+\frac{1 \pm \sqrt{3}\i}{2}=1. 
\]
This implies $b=(\pm 3+\sqrt{3}\i)/6$.
In the same way, from the conditions $t_{m_2m_4}(\rho_*)=t_{m_2m_5}(\rho_*)=1$, 
we obtain $f=j=(\pm 3+\sqrt{3}\i)/6$. 

Suppose $a=(1 + \sqrt{3}\i)/2$. Then we have 
\[
b=f=j=(3+\sqrt{3}\i)/6.
\] 
Furthermore, $t_{m_3m_4}(\rho_*)=t_{m_3m_5}(\rho_*)=t_{m_4m_5}(\rho_*)=1$ imply 
\[
ch+dg=-\frac{1}{3},\ cl+dk=-\frac{1}{3},\ gl+hk=-\frac{1}{3}.
\]
On the other hand, the determinant conditions yield 
\[
cd=-\frac{2}{3},\ gh=-\frac{2}{3},\ kl=-\frac{2}{3}.
\]
Combining these equations, we obtain 
\[
2g^2-cg+2c^2=0, \ 2k^2-ck+2c^2=0, \ 2k^2-gk+2g^2=0.
\]
The first two equations imply that $g$ and $k$ must be of the form $c(1 \pm \sqrt{15}\,\i)/4$.
Substituting these expressions into the third equation yields $c^2=0$, that is, $c=0$,
which contradicts the earlier condition $cd=-2/3 \neq 0$. 
The remaining case $a=(1 - \sqrt{3}\i)/2$ can be shown similarly. 
Therefore, the preimage $(h^*)^{-1}(\g_1)$ is empty, and hence the map $h^*$ is not surjective.

In contrast, for the ghost character $\g_2=(0,-1)$, 
we construct two representations $\rho_{\pm}: \pi_1(\Sigma_2T_{5,6}) \to \SL_2(\C)$ 
satisfying $h^*(\chi_{\rho_{\pm}})=\g_2$. Let $\i=\sqrt{-1}$. 
\[
\begin{array}{lll}
\rho_{\pm}(m_1m_2)=
\left(\begin{array}{cc} \i & 0\\
0 & -\i \end{array}\right)
&\Rightarrow& 
\tr(\rho_{\pm}(m_1m_2))=0,\\
\rho_{\pm}(m_1m_3)=
\left(\begin{array}{cc} -\dfrac{1}{2} & \dfrac{-2\i \pm \sqrt{5}}{4}\\
\dfrac{-2\i \mp \sqrt{5}}{3} & -\dfrac{1}{2} \end{array}\right)
&\Rightarrow& 
\tr(\rho_{\pm}(m_1m_3))=-1,\\
\rho_{\pm}(m_1m_4)=
\left(\begin{array}{cc} \dfrac{-1-\i}{2} & \dfrac{-\i \mp \sqrt{5}}{4}\\
\dfrac{-\i \pm \sqrt{5}}{3} & \dfrac{-1+\i}{2} \end{array}\right)
&\Rightarrow&
\tr(\rho_{\pm}(m_1m_4))=-1,\\
\rho_{\pm}(m_1m_5)=
\left(\begin{array}{cc} -\dfrac{\i}{2} & \dfrac{3\i}{4}\vspace*{0.1cm}\\
\i & \dfrac{\i}{2} \end{array}\right)
&\Rightarrow& \tr(\rho_{\pm}(m_1m_5))=0,  
\end{array}
\]
where the signs in the above descriptions are taken consistently. 
Note that the characters $\chi_{\rho_{+}}$ and $\chi_{\rho_{-}}$ represent 
distinct points in $(h^*)^{-1}(\g_2) \subset X(\Sigma_2T_{5,6})$, 
since $y_{234}(\chi_{\rho_+}) \neq y_{234}(\chi_{\rho_-})$. 
Thus, the map $h^*$ is not injective. This completes the proof.
\end{proof}

One can also verify that the preimage of the ghost character 
$(x_{12},x_{13})=(-2,1)$ under $h^*$ is empty. 
Moreover, the above representations $\rho_{\pm}$, together with 
Theorem 4.9 (2) in \cite{Nagasato5}, show that the map $\widehat{\Phi}$ 
is not surjective, as in the case of $T_{4,5}$. 
In other words, the fundamental group $\pi_1(\Sigma_2T_{5,6})$ admits 
an $\SL_2(\C)$-representation that cannot be realized, via the map $\widehat{\Phi}$ 
(in the sense of \eqref{rep_s2k}), from any trace-free $\SL_2(\C)$-representation 
of the knot group $G(T_{5,6})$. 


\section*{Acknowledgements}
The first author was partially supported, during the early stages of this research, 
by the JSPS Research Fellowships for Young Scientists, 
JSPS KAKENHI for Young Scientists (Start-up), 
MEXT KAKENHI for Young Scientists (B), 
and JSPS KAKENHI for Young Scientists (B). 
This work was also partially supported by JSPS KAKENHI (C), 
Grant Number 20K03619.


\end{document}